\documentstyle[a4wide,twoside,12pt]{article}
\def\k#1{\kern#1em}
\def\Ib#1{{I\kern-.25em#1}}
\def\Ibb#1{{I\kern-.23em#1}}
\def\vh#1{\vrule width.02em height#1ex depth0ex}
\def\vb#1{\vrule width.02em height1.47ex depth#1ex}
\def\vcg{\vrule width.02em height1.4ex depth-.05ex}
\def\NN{\Ibb N}
\def\ZZ{{\k{.26}\vh{0.5}\k{.04}\vb{-1}\k{-.34}Z}}
\newcommand{\DLZZ}{\settowidth{\longueur}{Z}Z\hspace{-0.9\longueur}Z}

\newcommand{\fer}[2]{\left\langle {{#1\atop #2 }} \right\rangle }

\newcommand{\ord}{\mbox{\rm ord}}
\begin{document}

\begin{titlepage}
\title{Maple Umbral Calculus Package}
\author{A. Bottreau\thanks{LaBRI,
Universit\'e de Bordeaux I,
33405 Talence Cedex, France,
{\tt bottreau@labri.u-bordeaux.fr}}
\and
A. Di Bucchianico\thanks{Eindhoven University of Technology,
Department of Mathematics and Computing Science,
P.O. Box 513,
5600 MB Eindhoven,
The Netherlands,
{\tt sandro@win.tue.nl}. Author supported by NATO CRG 930554.}
\and
D.\ E.\ Loeb\thanks{LaBRI,
Universit\'e de Bordeaux I,
33405 Talence Cedex,
France,
{\tt loeb@labri.u-bordeaux.fr}, URL: {\tt
http://www.labri.u-bordeaux.fr/$\sim$loeb/}. Author partially
supported by URA CNRS 1304, 
EC grant CHRX-CT93-0400, the PRC Maths-Info, and NATO CRG 930554.}}
%\date{\today}
\end{titlepage}
\maketitle

\begin{abstract}
Rota's Umbral Calculus uses sequences of Sheffer polynomials to count
certain combinatorial objects. We have developed a Maple package that
implements Rota's Umbral Calculus and some of its generalizations. A
Mathematica version of this package is being developed in parallel.

\begin{center}
{\bf R\'esum\'e}
\end{center}

Le calcul ombral de Rota utilise des suites de polyn\^omes de Sheffer
pour  
%\'enum\'er 
\'enum\'erer %%%NEW
certains objets combinatoires. Nous avons
d\'evelopp\'e une biblioth\`eque de fonctions Maple qui impl\'ementent cette
th\'eorie ainsi que ses g\'en\'eralisations. 
Une version Mathematica de ce package est d\'evelopp\'ee en
parall\`ele.\\[5mm] 


{\bf Hardware requirements}: A computer equipped with Maple {\bf V}.2
or {\bf V}.3.
\end{abstract}

\section{Introduction}

Umbral calculus is the study of the analogies between various
polynomial sequences and the powers sequence $x^{n}$. For example,
$x^{n}$ has many parallels with the lower factorial sequence
$(x)_{n}=x(x-1)\cdots (x-n+1)$:
\begin{itemize}
\item The forward difference operator $\Delta:p(x)\mapsto p(x+1)-p(x)$ plays
a role with respect to $(x)_{n}$ analogous to that played by the
derivative $d$ with respect to $x^{n}$. 
\item Taylor's theorem is analogous
to Newton's theorem.
\item The binomial theorem for $(x+y)^{n}$ is replaced by Vandermonde's
identity for $(x+y)_{n}$.
\end{itemize}

Although Umbral Calculus dates back to the 18th century, it 
was only put on a rigorous foundation by Gian-Carlo Rota and his
collaborators \cite{MuR,RKO}  in the 1970's. We now characterize each
polynomial sequence 
under study by one or more polynomial operators (usually
shift-invariant\footnote{Shift-invariant operators commute with the
shift operator $E^{a} \colon p(x) \rightarrow p(x+a)$.}) 
associated with it. 
The duality between operators and 
polynomials is the key tool to deriving umbral calculus results.

Umbral Calculus has many applications in enumerative combinatorics. 
The powers $x^{n}$ counts all functions from an $n$-element set to an
$x$-element set, while the lower factorials $(x)_{n}$ count
injections. 
Similarly, given any species of combinatorial structures (or
quasi-species), let $p_{n}(x)$ be the number of functions from an
$n$-element set to an $x$-element set enriched by this species. A
function is enriched by
associating a (weighted) structure with each of its fibers. 
The resulting sequence of polynomials 
$\left(p_{n}\right)_{n\in \NN} $ is said to be of 
binomial type since it obeys the ``binomial'' identity
$$ p_{n}(x+y)=\sum _{k=0}^{n} {n\choose k }p_{k}(x) p_{n-k}(y). $$
For example,
given the species of rooted forests, the enriched functions are called
persistent functions and are enumerated by the Abel polynomials
$A_{n}(x)=x(x+n)^{n-1}$.
Other applications include lattice path counting
\cite{Nie11,Raz1,Wat3}. 

Our Maple package provides a number of different tools by which
to enter operators. These operators can then be manipulated in many
different ways. In particular, the polynomial sequences associated
with them can be explicitly calculated.

This package has already aided us in our research \cite{pers}; 
we hope that it will help you too.

We expect to release a Mathematica version of this package in the near
future.  

\section{Polynomial Operators}

Polynomial operators (shift-invariant or not) can be specified in
several convenient manners:
\begin{description}
\item[Explicitly by their action on polynomials.] For example, the
shift operator is defined
{\tt < subs(x = x+a, p) | p >}
using the ``angle-bracket'' notation for functional operators.
(See {\tt ?operators[functional]} and {\tt ?unapply} for details.)
Similarly, the Bernoulli operator
$p(x)\mapsto \int_{x}^{x+1}p(t)dt$
is defined {\tt <int(subs(x=t, p), t=x..(x+1)) | p | t>}~.
\item[As an analytic function of the derivative.]
By the expansion theorem 
(see \cite[Theorem~2, p.~185]{MuR} or \cite[Theorem~2, p.~691]{RKO}), %%% NEW
any shift-in\-var\-i\-ant operator can be
expanded as a formal power series in {\tt d}  where {\tt d} is a special 
reserved symbol representing the derivative. 
For example, the shift operator is defined 
{\tt exp(a*d)}.
\item[Abstractly as an unspecified function of {\tt d}.] For example, {\tt
f(d)}  or {\tt f(d,x)} in the case of a non-shift-invariant operator
\cite{KuM}. 
\item [Using the {\tt powseries} package.] If the
coefficients of the 
formal power series given by the expansion theorem are all known, then
use {\tt powcreate}. For example, 
{\tt powseries [powcreate] (f(n) = a\^{}n/n!);}. 
\item [As a series.] If only finitely many terms are known, then use
{\tt series}. For example, 
{\tt series(1 + a*d + a\^{}2*d\^{}2/2 + c*d\^{}3, d, 3);}.
\end{description}
See {\tt ?linear} for more information.

Operators can be converted easily from one form to another with {\tt
convert}. A delta operator {\tt Q} is a shift-invariant operator such that
{\tt Q (x)} is a non-zero constant. An abstract operator is assumed to be
invertible unless indicated otherwise. (For example, {\tt Q :=
d*f(d);} or {\tt Q := f(d); f(0) := 0;}.)

Using the generalized expansion theorem, shift-invariant operators can
be expanded with {\tt operatorExpansion}  as a formal power series
in an arbitrary delta operator. 
\begin{verbatim}
> oe(d,delta,3);
                                            2                 3
               exp(d) - 1 - 1/2 (exp(d) - 1)  + O((exp(d) - 1) )
\end{verbatim}
Such operator expansions are practical for numerical calculations. 
For example, expanding the Bernoulli operator  into powers
of the forward difference operator {\tt delta} yields the classical
Newton-Cotes formula of numerical 
integration \cite[p.\ 186]{MuR}. 

Our package also allows the expansion of linear operators which are
not shift-invariant. Such expansions \cite{KuM} express the operator
as a formal power series in {\tt d} whose coefficients are polynomials
in {\tt x}.
For example, if {\tt Q} is the operator $ {\tt Q} \colon p(x)\mapsto
\int_{0}^{x}p(t)dt$, then 
the expansion {\tt convert(Q,function,5,x)} or {\tt 
convert(Q,powseries,x)} of {\tt Q} in terms of multiplication by {\tt
x} and the derivative {\tt d} gives an elementary proof of Bourbaki's
method of  asymptotic integration \cite[Sections 3.5 and 3.6]{Bour}.

Operators can be applied to polynomials with {\tt dp} by specifying
the free variable of the polynomial:
\begin{verbatim}
> dp(delta,x^3,x);
                                    2
                                 3 x  + 3 x + 1
\end{verbatim}
In case the degree of the polynomial is not explicitly given, the
program will calculate the {\tt Order} most significant terms of
the answer. ({\tt Order} is a system variable whose default value is
6.) 
\begin{verbatim}
> dp(delta,x^n,x);
     (n-1)                  (n - 2)                      (n-3)
  n x        + 1/2 n (n-1) x        + 1/6 n (n-1) (n-2) x
                                   (n-4)
       + 1/24 n (n-1) (n-2) (n-3) x
                                          (n-5)      (n-6)
       + 1/120 n (n-1) (n-2) (n-3) (n-4) x      + O(x     )
\end{verbatim}

\section{Polynomial Sequences}

Given the necessary operators, the program can calculate
polynomials of binomial type ({\tt bfo}), Sheffer 
sequences ({\tt sfo}), Steffensen sequences ({\tt steff}), and
cross-sequences ({\tt cseq})
(see \cite[Sections~5 an 8]{RKO}). %%% NEW
For example, the sequence of binomial type
$ {\tt bfo(p(d),x,n)}$ associated with a delta operator {\tt
p(d)} is defined by the conditions 
$$\begin{array}{rcll}
{\tt dp(p(d),bfo(p(d),x,n))} & {\tt =} & {\tt n * bfo(p(d),x,n)} &
\mbox{for {\tt n$>$0}}\\
{\tt bfo(p(d),0,n)} & {\tt =} & {\tt 0} &
\mbox{for {\tt n$>$0}}\\
{\tt bfo(p(d),x,0)} & {\tt =} & {\tt 1},
\end{array}$$
or equivalently by its exponential generating
function $\exp(q(x)t)$ where
{\tt q} is the compositional inverse of {\tt p}.  
Note that {\tt q(t)} is the generating
function of the associated species.

For example, the lower
factorial $(x)_{n}$ is the basic sequence of binomial type for the
forward difference operator $ \Delta$.
\begin{verbatim}
> factor(bfo(delta,x,4));
                           x (x - 1) (x - 2) (x - 3)
\end{verbatim}
If the degree is not explicitly given, then only the most significant
terms will be computed.

Several functions in the package do further operations on polynomial
sequences.  Arbitrary polynomials can be expressed in terms of such
sequences ({\tt polynomialExpansion, shefferExpansion,
basicExpansion}).  For example, \begin{verbatim}> p := randpoly(x);
                          5       4       3       2
                 p := 79 x  + 56 x  + 49 x  + 63 x  + 57 x - 59
> be(delta,p,x);
         - 59 bfo(exp(d) - 1, x, 0) + 304 bfo(exp(d) - 1, x, 1)
              + 1787 bfo(exp(d) - 1, x, 2) + 2360 bfo(exp(d) - 1, x, 3)
              + 846 bfo(exp(d) - 1, x, 4) + 79 bfo(exp(d) - 1, x, 5)
> ";
                       5       4       3       2
                   79 x  + 56 x  + 49 x  + 63 x  + 57 x - 59
\end{verbatim}
Connection constants can be determined between arbitrary polynomial
sequences. For example, the Stirling numbers are given by 
{\tt cc(topseq(powerx, 5, x), topseq(lower, 5, x), x)} where {\tt
powerx(n,x)} is $x^{n}$ and {\tt lower(n,x)} is $(x)_{n}$. 
Other features include umbral composition ({\tt uc}), and umbral
inversion ({\tt ui}).

\section{Generalizations}

Several authors (eg. \cite{Rom9,Vis1}) have generalized the umbral
calculus by considering not only sequences of binomial type with
generating function $\exp(g(x)t)$ but also those whose generating
function is $\Phi(g(x)t)$ where $\Phi(t)=\sum _{n=0}^{\infty
}t^{n}/[n]!$ and $[n]!$ denotes the generalized factorial 
$[n]!={\tt a(1)a(2)\cdots a(n)}$. Most of the functions in the umbral calculus package
allow an optional argument {\tt a} which is either left undefined, or
defines the coefficients used
by the ``generalized derivative.'' Thus, 
\begin{verbatim}
>  dp(d, x^3, x, proc(n) 1 end);
                              2
                             x 
\end{verbatim}
The following possible choices for {\tt a} are predefined.
\begin{center}
\begin{tabular}{|rl|c|r@{ {\tt =} }l|}
\hline
\multicolumn{2}{|c|}{Umbral Calculus}&{\tt dp(d,p,x,a)}&
\multicolumn{2}{|c|}{{\tt a}}\\[2mm] 
\hline
Classical&\cite{MuR,RKO}&{$\displaystyle \frac{dp(x)}{dx}$}&{\tt
classical(n)}&{\tt n} \\[6mm] 
$q$-umbral calculus&\cite{Rom9,Rom10}&{$\displaystyle
\frac{p(qx)-p(x)}{(q-1)x}$}&{\tt 
gaussian(n)}&{\tt (q\^{}n-1)/(q-1)}\\[6mm]
Divided Difference&\cite{HiR,VeS3}&{$\displaystyle \frac{p(x)-p(0)}{x}$}&{\tt
divided(n)}&{\tt 1}\\[6mm]
Hyperbolic&\cite{pers}&{$\displaystyle
\left(\frac{d}{d\sqrt{x}}\right)^{2} p(x)$}&{\tt   
hyperbolic(n)}&{\tt 2*n*(2*n-1)}\\ %%% NEW I put braces around \sqrt{x}
\hline
\end{tabular}
\end{center} 
See {\tt ?genderiv} for details.

Generalizations of the umbral calculus to several variables 
\cite{Rom3,Wat3} are supported. Most functions included
in the package have an alternate syntax for use in multivariate umbral
calculi. In particular, {\tt d[i]} represents the partial derivative with
respect to the $i$th variable. Instead of a single delta operator, a
collection of operators are required to define a sequence of binomial
type. This generalization is completely compatible with the above
generalization. See {\tt ?multilinear} and {\tt ?moe} for details.

For further instructions consult the on-line help and examples provided in
the package. For help, type {\tt ?}{\em key-word}. An index of
key-words is available via {\tt ?umbral}.

See \cite{survey} for an extensive survey and  bibliography of the
umbral calculus. 

%\bibliography{umbral,par}

\end{document}